\documentclass[12pt,a4paper]{amsart}
\usepackage[latin1]{inputenc}
\usepackage{amsmath}
\usepackage{amsfonts}
\usepackage{amssymb}
\usepackage{amsmath,bbm,amssymb,amsxtra}
\usepackage{color}
\usepackage{enumerate}
\allowdisplaybreaks
\usepackage{epsfig}
\usepackage{ae}
\newcommand{\refeq}[1]{(\ref{#1})}

\theoremstyle{definition}

\newtheorem{lemma}{Lemma}[section]
\newtheorem{proposition}[lemma]{Proposition}

\newtheorem{theorem}[lemma]{Theorem}

\newtheorem{remark}[lemma]{Remark}
\newtheorem{definition}[lemma]{Definition}

\newcommand{\prop}[1]{\begin{proposition}\label{#1}
\sl }
\newcommand{\eprop}{\end{proposition}}
\newcommand{\thm}[1]{\begin{theorem}\label{#1}
\sl }
\newcommand{\ethm}{\end{theorem}}
\newcommand{\lem}[1]{\begin{lemma}\label{#1}
\sl }
\newcommand{\elem}{\end{lemma}}

\newcommand{\defin}[1]{\begin{definition}\label{#1}
\sl }
\newcommand{\edefin}{\end{definition}}

\newcommand{\beqno}{\begin{eqnarray*}}
\newcommand{\eeqno}{\end{eqnarray*}}
\newcommand{\beqla}[1] {\begin {eqnarray}\label{#1}}
\def\eeq {\end {eqnarray}}
\newcommand{\beq}{\begin {eqnarray}}

\newcommand{\real}{{\mathbb R}}

\newcommand{\integer}{{\mathbb Z}}
\newcommand{\nanu}{{\mathbb N}}
\newcommand{\complex}{{\mathbb C}}

\newcommand{\torus}{{\mathbb T}}

\author[Drasin]{David Drasin}
\address{Department of Mathematics, Purdue University, West Lafayette,
IN 47907, USA}
\email{drasin@math.purdue.edu}
\title[Frequently hypercyclic entire functions]{Optimal growth of frequently hypercyclic entire functions} 
\date{September 2, 2011}
\author[Saksman]{Eero Saksman}
\address{Department of Mathematics and Statistics, University of
Helsinki, PO~Box~68, FI-00014 Helsinki, Finland}
\email{eero.saksman@helsinki.fi}

\subjclass{Primary 47A16, secondary  30D15, 47B38}

\keywords{Frequently  hypercyclic operator, differentiation operator, rate of growth, entire functions}

\thanks{The first author thanks the University of Helsinki for its hospitality and stimulation provided in 2005--6 during his sabbatical, when the work for part (i) of the main result was obtained. It was during this period that both authors attended a workshop in Functional Analysis sponsored by the University of M\'alaga  \cite{BG-E1} where we learned of this interesting problem. The second  author was
supported by the Finnish CoE in Analysis and Dynamics Research,
and by the Academy of Finland, projects 
113826 \& 118765}
\begin{document}

\centerline{\it Dedicated to Prof. Olli Martio on the occasion of his 70'th birthday}

\smallskip

\begin{abstract} We solve a problem posed by A. Bonilla and K.-G. Grosse-Erdmann \cite{BG-E1} by constructing an entire function $f$ that is frequently hypercyclic with respect to the differentiation operator, and satisfies
$M_f(r)\leq\displaystyle  ce^r r^{-1/4}$, where  $c>0$  be chosen arbirarily small. The obtained growth rate is sharp. We also obtain optimal results for the growth when measured in terms of average $L^p$-norms. Among other things, the proof applies Rudin-Shapiro polynomials and heat kernel estimates. 
\end{abstract}
\maketitle

\section{introduction}\label{se:intro}

A linear operator $T$ on a separable topological vector space $E$.
is {\sl hypercyclic} if there exists  $x\in E$ such that the set of  iterates $\{ T^nx : n\geq 1\}$ is dense in $E.$
In this situation $x$ is sometimes called an {\sl universal} element.
We refer to \cite{G-E} for basic facts on hypercyclicity.

Recently there has been  interest on a related, more stringent notion. The operator $T$ (and likewise  the element $x\in E$) is called {\it frequently hypercyclic} if $T^nx$ visits any given neighbourhood 
with a relatively constant rate. More precisely, given any open set $U\subset E$ one asks that the set
$$
A=\{n\geq 1\; :\; T^nx\in U\}
$$
has positive density, i.e.
$$
\liminf_{n\to\infty}\frac{1}{n}\#(A\cap\{ 1,\ldots ,n\}) >0.
$$
This notion was introduced by Bayart and Grivaux \cite{BG1}, \cite{BG2} and has been studied in many papers devoted to operators in Hilbert, Banach, or general topological vector spaces. We  refer to \cite{G} and \cite{BBG-E} and especially for the references therein for more information on the known results.

Classical examples of hypercyclic operators are  the translation and differentiation operators in the space $\mathcal E$  of entire functions on the complex plane $\complex$, equipped with the standard compact-open topology. We shall consider only the differentiation operator $D:
\mathcal E\to \mathcal E$, where $Df(z):=f'(z)$ (see \cite{BBG-E} for results on the translation operator). Let us recall that hypercyclicity of $D$ is a classical result due to MacLane  \cite{MacL}.

In this note we study the following problem: how slowly can a D-frequently hypercyclic entire function grow? This question was raised by Bonillla and Grosse-Erdmann  \cite{BG-E1}, and at the same time they gave concrete estimates for the minimal growth of such a function (it was proven in  \cite{BG-E2} that indeed $D$  is frequently hypercyclic).  In \cite{BG-E0} the same authors  generalized a well-known result of Godefroy and Shapiro \cite{GS} by showing that every operator on $\mathcal E$ which commutes with $D$, and  is  not a multiple of the identity, is frequently hypercyclic. A couple of years later  Blasco, Bonilla and  Grosse-Erdmann \cite{BBG-E}
improved on  the earlier results  and 
 showed that a  $D$-frequently hypercyclic entire function $f$ satisfies
\beqla{necessary}
\liminf_{r\to\infty} \frac{M_f(r)}{r^{1/4}}>0,
\eeq
where  $M_f(r):=\sup_\theta|f(re^{i\theta})|$.
Moreover, given any function $\phi:\real\to[1,\infty )$ with
$\lim_{r\to\infty}\phi (r)=\infty $ they proved  the existence of a  $D$-frequently hypercyclic entire function $f$ with
\beqla{sufficient}
M_f(r)\leq e^r\phi (r) \qquad {\rm for}\;\; r\geq 1.
\eeq

The paper \cite{BBG-E} also considered growth in terms of the average $L^p$-norms. Thus,   for $p\in [1,\infty )$  and $r>0$ let
$$
M_{f,p}(r):= \big( \frac{1}{2\pi} \int_0^{2\pi}
 |f(re^{i\theta})|^p\, d\theta\big)^{1/p},
$$
and augment this with $M_{f,\infty}(r) =M_f(r).$
In \cite{BBG-E}, the authors also  showed  that, given $a>0$, an inequality
\beqla{inequality}
M_{f,p}(r)\leq C_a e^r r^{-a}\qquad \nonumber
\eeq
could only be valid if
$a\leq \max (1/2p, 1/4)$. Moreover, if $C$ is replaced by a factor $\phi(r)\uparrow\infty,$
examples with the growth rate   $a=\min (1/2p, 1/4)$ were constructed. 

Quite recently Bonet and Bonilla \cite{Bo2} constructed almost optimal examples in the range $p\in [1,2)$ by showing that $a=1/2p$ can be achieved, again requiring the factor $\phi(r)\uparrow\infty$. Their result is  sharp in the special case $p=1$.

Our main result  determines the optimal growth rate of entire $D$-frequently hypercyclic functions. It turns out that the sharp result actually corresponds to the lowest possible rates consistent
with \cite{BBG-E}  for all $p$.
Both \cite{BBG-E} and  \cite{Bo2} employ fairly sophisticated  tools from the general theory of frequent hypercyclicity. In contrast, our construction relies only on basic complex analysis, elementary  heat kernel estimates and on   two classes of classical
polynomials whose properties reflect the role that $p$ plays in 
these results.    Thus for large $p$, our construction patches together Rudin-Shapiro polynomials $p_m$, having coefficients $\pm 1$, but
whose supnorm is minimal; for $1\leq p\leq 2$ an analogous role is played by the de la Vall\'ee-Poussin    polynomials $p_m^*$.

We present our main result in three parts, with the most interesting case $p=\infty$ meriting its
own statement. 
\thm{th:main}\quad 
{\bf (i)}\quad For any $c>0$ there is an entire frequently hypercyclic function $f$ such that
$$
M_f(r)\leq c\frac{e^{r}}{r^{1/4}}\qquad {\rm for \; all }\;\, r>0.
$$
This estimate is optimal: every such function satisfies $\limsup_{r\to\infty} r^{1/4}e^{-r}M_f(r)>0.$
\smallskip

\noindent{\bf (ii)}\quad More generally, 
given  $c>0$ and $p\in (1,\infty]$ there is an entire $D$-frequently hypercyclic function $f$ with
$$
M_{f,p}(r)\leq c\frac{e^{r}}{r^{a(p)}}\qquad {\rm for \; all }\;\, r>0,
$$
where $a(p)=1/4$ for $p\in [2,\infty]$ and $a(p)=1/(2p)$ for $p\in (1,2].$
This estimate is optimal: every such function satisfies $\limsup_{r\to\infty} r^{a(p)}e^{-r}M_{f,p}(r)>0.$

\noindent{\bf (iii)}\quad Given $\phi(r)\uparrow\infty$, there is an entire $D$-frequently hypercyclic function $f$
with
$$
M_{f,1}(r)\leq \phi(r)\frac{e^{r}}{r^{1/2}}\qquad {\rm for \; all}\;\,r >0.
$$
This estimate is optimal: every such function satisfies $\limsup_{r\to\infty} r^{1/2}e^{-r}M_{f,1}(r)=\infty.$
\ethm
\noindent
The sharp conclusion (i)  provides the converse of \refeq{necessary}, yielding a considerable  strengthening of the best previously known growth \refeq{sufficient}. In turn, (ii) 
sharpens the main result of \cite{Bo2} by removing the unnecessary increasing factor $\phi (r)\uparrow \infty$. Finally, assertion (iii) is already due to Bonet and Bonilla   \cite[Corollary 2.4]{Bo2}.

\section{The construction of the $D$-frequently-hypercyclic function $f$}
For any given polynomial $q$ with Taylor series
$$
q(z)=\sum_{j=0}^d\frac{q_jz^j}{j!},\qquad d={\rm deg}\, (q),
$$
set
$$
\widetilde{q}(z):= \sum_{j=0}^d{q_jz^j}\qquad{\rm and}\quad \|\widetilde{q}\|_{\ell^1}:=\sum_{j=0}^d|q_j|. 
$$

For our purposes it will be useful to divide  the Taylor series  $
\displaystyle f(z)=\sum_{k=0}^\infty \frac {a_k}{k!}z^k,
$ of a given  entire function $f$,
 into blocks of specified size. It turns out that the correct partition is the decomposition
$$
f(z)=\sum_{n=0}^\infty P_nf,
$$
where 
\beqla{blocks}
P_nf(z)=\sum_{k=n^2}^{(n+1)^2-1} \frac {a_k}{k!}z^k\qquad(n\geq 0).
\eeq
In fact, the precise  size of the blocks $P_k$ is one of the key  ingredients of our argument. The motivation for the choice \refeq{blocks} will be discussed later on in Remark \ref{re:reason}.

\subsection{Long polynomials with controlled norm} 
Explicit polynomials  having  a fixed proposition of coefficients 1 and with small $L^p$ -norm (up to order
of magnitude) will be  needed for construction of the blocks $P_kf$. It is here  where the cases $1\leq p\leq 2$ and $p\geq 2$ bifurcate. 
The first part of
the next lemma records the beautiful  result of Rudin and Shapiro \cite{R} which produces polynomials $\{p_m\}$
of each degree $m$ having coefficients $\pm1$ and with optimal growth of sup norm (see also \cite[\S 6, Chapter 5 ]{K1}, \cite{K0} and \cite{BB}  for further results on unimodular polynomials). That this growth is indeed optimal  follows immediately by the Parseval formula.
For $p\in [1,2)$,
the $\{p_m\}$ are replaced by the de la Vall\' ee-Poussin polynomials $p^*_m$.  Below $p'$ stands for the  exponent conjugate to $p$.

\lem{le:rudin-shapiro}
{\bf (i)}\quad For each $m\geq 1$ there is a trigonometric polynomial $p_m$ 
$$
p_m=\sum_{k=0}^{m-1} b_{m,k}e^{ik\theta},
$$
where $b_{m,k}= \pm 1$ for all $0\leq k\leq m-1$,  with at least half of the coefficients of $f$ being $+1$, and with
$$
\|p_m\|_{p}\leq 5\sqrt{m}\quad {\rm for} \;\; p\in [2,\infty].
$$

\noindent {\bf (ii)}\quad  Corresponding to each  $m\geq 1$ is a polynomial
$$
p^*_m=\sum_{k=0}^{m-1} b^*_{m,k}e^{ik\theta},
$$
where $|b_{m,k}|\leq 1$ for all $0\leq k\leq m-1$, and
with  at least $\lfloor m/4\rfloor$  coefficients being $+1$, and with
$$
\|p^*_m\|_{p}\leq 3m^{(1/p')} \quad {\rm for } \;\; p\in [1,2].
$$
\elem
\begin{proof} In case  $p=\infty$ assertion (i)  just records the main result of \cite{R}. Since $\| p_m\|_2=\sqrt{m},$ the claim for exponents $p\in (2,\infty)$ is immediate. In turn,  for assertion (ii) we may assume that $m\geq 4$. Set $k=\lfloor m/4\rfloor$ 
and choose let $p^*_m(e^{i\theta})=e^{2ki\theta}(2F_{2k}(e^{i\theta})-F_k(e^{i\theta}))$, where $F_k$ is the $k$:th Fejer kernel. In other words, $p^*_m$ corresponds to a shifted de la Valle Poussin kernel \cite[p.16]{K1}. Then $\|p^*_m\|_1\leq 3$ and the statement for  $p\in (1,2]$ is obtained by interpolating (i.e. applying H\"older's inequality)  with the obvious estimate $\|p^*_m\|_2\leq \sqrt{m}.$
\end{proof}

\subsection{Explicit formula for $f$}\label{f} Let $0<c<1$ be given, let $p\in (1,\infty ]$.  Denote by $\mathcal{P}$  the (countable) set of all polynomials with rational
coefficients, and consider pairs $(q,\ell)$ with $q\in\mathcal{P}$ and $\ell\in\mathbb{N}$ with
$\ell\geq \|\widetilde{q}\|_{\ell^1},$ exhibited in a single sequence $(q_k,\ell_k)\,{k\geq 1}$. Let us record the important fact
\beqla{imp}
 \|\widetilde{q_k}\|_{\ell^1}\leq \ell_k\quad {\rm for \; every}\;\; k\geq 1.
\eeq

Next  partition  the even integers in $\nanu$ into countably many infinite disjoint arithmetic
sequences by setting $2\nanu =\cup_{k\geq 1}{\mathcal A}_k$ where
$$
{\mathcal A}_k=\{ 2^k(2j-1) \; :\; j\in\nanu\} .
$$
Next, for any $k\geq 1$ denote by $\alpha_k$ the integer
\beqla{alpha}
\alpha_k&:=& 1+ \big\lfloor\max \Big( ( 10^{10}\ell_k/c)^{\max (2,p')}, 2d_k+8\ell_k\Big )\big\rfloor ,\label{c_k}
\eeq
with $d_k$ the degree of $q_k$.

We have noted that
$$
f=:\sum_{j=1}^\infty \frac{a_jz^j}{j!}:=\sum_{n=1}^\infty P_nf
$$
is uniquely determined by the blocks $P_nf.$ 
First, consider the case $p\geq 2$. Set first $a_0=0$, or, in other words $P_0f=0.$ Fix
$n\geq 1$. If $n$ is odd, take $P_nf=0$. When $n$ is even and  $n\in {\mathcal A}_k$ we use the Rudin-Shapiro polynomials to write
\beqla{piece}
\widetilde{P_nf}:= \left\{ 
\begin{array}{ll}
0& {\rm if}\quad n<10\alpha_k,\\
z^{n^2}p_{\lfloor n/\alpha_k\rfloor}(z^{\alpha_k})\widetilde{q_k}(z) & {\rm otherwise}.
\end{array}
\right.
\eeq

For $1<p<2,$  one just modifies the second line in \refeq{piece}  by setting for
even $n$, $n\in\mathcal{A}_k, n\geq 10\alpha_k$
\beqla{piece2}
\widetilde{P_nf} = z^{n^2}p^*_{[n/\alpha_k]} (z^{\alpha_k})\widetilde{q_k}(z).
\eeq 

We shall treat the case $p=1$ in Remark \ref{p=1} below. In what follows we shall show that the function $f$ we just defined satisfies the assertions (i) and (ii) of Theorem \ref{th:main}.

\section{The proof}\label{se:proof}
The proof of Theorem \ref{th:main}  will be based on two auxiliary results, Propositions \ref{glue} and \ref{heat} below. We begin by recording  two simple  auxiliary observations.

\lem{sum}\label{le1}
 Let $a\in [0,1].$ Then for any  $m \geq 1$,
$$
\sum_{n=1}^\infty  e^{n^2}n^{-2a}\left( \frac{m^2}{n^2}\right)^{n^2}
\leq 10\, e^{m^2}m^{-2a}.
$$
\elem
\begin{proof}
Write the elementary inequality $1+\log x\leq x$ ($0< x\leq \infty$) as
$1+\log (x^2)\leq x^2-(1-x)^2.$ By making the substitution $x=m/n $ and
multiplying both sides by $n^2$ we obtain
$$
n^2+n^2\log (m^2/n^2)\leq m^2-(n-m)^2.
$$
Thus $e^{n^2}(m^2/n^2)^{n^2}\leq e^{m^2}e^{-(n-m)^2}$, and hence the sum has the upper bound
\beqla{eq111}
&&
 e^{m^2}m^{-2a}\Bigl( \sum_{n=1}^{m}(m/n)^{2a}e^{-(n-m)^2}+\sum_{n=m+1}^\infty e^{-(n-m)^2}\Bigr)
 \nonumber\\&&\leq
 e^{m^2}m^{-2a}\Bigl(2\sum_{j=0}^\infty (j+1)^2e^{-j^2}\Bigr) \leq 10\,  e^{m^2}m^{-2a},\nonumber
\eeq
where we have used the simple inequality $(m/n)^{2a}\leq ((m-n)+1)^{2a}\leq ((m-n)+1)^{2}$ when  $n\leq m.$ 
\end{proof}

\lem{loglinear} Assume that  $0<x_0<x_1$  and the function $u:[x_0,x_1]\to \real$ is of the form 
$$
u(x)=a\log (x)+b-(cx+d)
$$
 with $a>0$ and $u(x_0)=u(x_1)=0$. Then $u(x)\leq a(x_1/x_0 -1)^2/8$ for
all $x\in [x_0,x_1].$
\elem
\begin{proof}
Simply observe that $u(x)-(x -x_0)(x_1-x)ax_0^{-2}/2\leq 0$ as  the left hand side is convex and vanishes at the endpoints.
\end{proof}

The following proposition contains an important underlying principle for bounding the growth of a given entire function $g$: correct  size for each $P_ng$ on  both $\{ |z| =n^2\}$ and on  $\{ |z| =(n+1)^2\}$ is enough to quarantee the desired  growth for $g$.

\prop{glue} Let $p\in (1,\infty ]$ and $a\in [0,1].$  Assume that there is a constant $b>0$ such that  for each $n\geq 1$ the blocks
of  a given entire function $g$ with $g(0)=$ satisfy
\beqla{eq1}
M_{P_ng,p}(n^2)\leq b e^{n^2}n^{-2a}.
\eeq
and
\beqla{eq2}
M_{P_ng,p}((n+1)^2)\leq b e^{(n+1)^2}(n+1)^{-2a}.
\eeq
Then $g$ itself satisfies 
$$
M_{g,p}(r)\leq  10^3b e^{r}r^{-2a}\qquad (r>0).
$$
\eprop

\begin{proof} 
We start by recalling  \cite[p.76]{HK} that for any entire function $g$, the $p$-means $\log(M_{g,p}(r))$, $p>0$ are increasing, convex functions of $\log r$ (alternatively, in our range $p\in [1,\infty  ] $ this is just  Hadamard's three circles theorem combined with a  simple duality argument). Especially, one may apply the maximum principle for the $p$-means.
When $h=P_ng$, we have
that $z^{-n^2}h(z)$ is holomorphic in $\{ |z|\leq n^2\}$ and that $z^{-(n+1)^2}h(z)$ is holomorphic in
$\overline{\complex}\setminus \{ |z|<( n+1)^2\}$. Hence (\ref{eq1}) and (\ref{eq2}) imply by the maximum principle that
\beqla{eq3}
M_{P_ng,p}(m^2)\leq b e^{n^2}n^{-2a}\left(\frac{m^2}{n^2}\right)^{n^2}\qquad
{\rm for} \;\; 1\leq m\leq n.\nonumber
\eeq
and
\beqla{eq4}
M_{P_ng,p}(m^2)\leq b e^{(n+1)^2}(n+1)^{-2a}\left(\frac{m^2}{(n+1)^2}\right)^{(n+1)^2}\qquad
{\rm for} \;\; 1\leq n<m. \nonumber
\eeq

By summing over $n$, observing that $P_0g=0,$ and and invoking Lemma \ref{sum}
 we obtain  the desired estimate when $r=m^2$, with $m$ a positive integer
\beqla{eq5}
M_{g,p}(m^2)&\leq &\sum_{n=m}^\infty b e^{n^2}n^{-2a}\left(\frac{m^2}{n^2}\right)^{n^2}+ \sum_{n=1}^{m-1}  b e^{(n+1)^2}(n+1)^{-2a}\left(\frac{m^2}{(n+1)^2}\right)^{(n+1)^2}\nonumber\\
&\leq & 11\, be^{m^2}m^{-2a}.\nonumber
\eeq

The log-convexity  of the $p$-means together with Lemma \ref{loglinear} allow this estimate to be interpolated inside each interval $ (m^2 ,(m+1)^2)$.
Observe first that  the effect of the term $m^{-2a}$ can be ignored  since  $ \log (r^{-2a})$ is  a linear function of
$\log r.$ Thus assume that $a=0$ and $m\geq 1.$ Denote $r_0=m^2$ and $r_1=(m+1)^2$. We obtain
for $r\in (r_0, r_1)$ that
$$
M_{g,p}(r)\leq  \log 11 +\frac{\log r_1-\log r}{\log r_1-\log r_0}r_0+\frac{\log r-\log r_0}{\log r_1-\log r_0}r_1:=C
+\frac{r_1-r_0}{\log(r_1/r_0)}\log r.
$$
If one replaces  the right hand side  by a function that is linear in $r$ with the same values
at the endpoints $r_0,r_1,$ Lemma \ref{loglinear} yields that the error so induced is less than 
$$
\left(\frac{r_1-r_0}{\log (r_1/r_0)}\right)(r_1/r_0-1)^2/8\leq \frac{3(r_1-r_0)^2}{8r_0}\leq (3/2)^3,
$$
where we have also applied  the  inequality $\log (r_1/r_0)\geq (r_1/r_0 -1)/3$.The stated result  follows by observing that $11e^{(3/2)^3}\leq 10^3.$ 

Finally, the claim for values $r\in (0,1]$ follows immediately from the maximum principle as we already know that $M_{g,p}(1)\leq 11b$ and since $11e\leq 1000$.
\end{proof}

The relation between the degree-$2n$ polynomial determining $P_ng$ and the growth of $P_ng$ is characterized in the following proposition, where the main trick is the surprising appearance of the heat kernel.
\prop{heat} Assume that  $n\geq 1,$ $p\in [1,\infty ]$ and  $g$ is an arbitrary entire function
$\displaystyle g(z)=\sum_{k=0}^\infty \frac {b_k}{k!}z^k.$ Then, if
\beqla{eq219}
B:=\|\sum_{k=0}^{2n}b_{n^2+k}e^{ik\theta}\|_{L^p(\torus )}.\nonumber
\eeq
one has
\beqla{eq20}
M_{P_ng,p}(n^2)\leq  20\, B e^{n^2}n^{-1}
\eeq
and
\beqla{eq21}
M_{P_ng,p}(n+1)^2)\leq 20\, B e^{(n+1)^2}(n+1)^{-1}.
\eeq

\eprop

\begin{proof}
We consider first \refeq{eq20}, and   may  assume that $n\geq 2.$ Observe that one may write
\beqla{eq22}
|P_ng(n^2e^{i\theta} )|=\frac{(n^2)^{n^2}}{(n^2)!}\Big|
\sum_{k=0}^{2n}\lambda_{n,k}b_{n^2+k}e^{ik\theta}
\Big|, 
\eeq
where
\beqla{eq23}
\lambda_{n,k}:= \left(\frac{n^2}{n^2}\right)\left(\frac{n^2}{n^2+1}\right)\left(\frac{n^2}{n^2+2}\right)\cdots
\left(\frac{n^2}{n^2+k}\right).
\eeq
Stirling's formula yields the rate of growth:
\beqla{eq24}
\frac{(n^2)^{n^2}}{(n^2)!}\leq \frac{e^{n^2}}{\sqrt{2\pi} n}.\nonumber
\eeq
Hence, we need to verify that the Fourier multiplier operator
\beqla{eq25}
\sum_{k=0}^{2n-1}c_{k}e^{ik\theta}\mapsto  \sum_{k=0}^{2n-1}\lambda_{n,k}c_{k}e^{ik\theta},\nonumber
\eeq
acting on the $L^p$-space of trigonometric polynomials of degree $2n$, has norm bounded
independent of $n$ and $p.$ 

For that end we first compute (recall that $k\leq 2n$)
\beqla{eq26}
\ \ \ -\log(\lambda_{n,k} )&=&\sum_{j=1}^k\log \Big(\frac{n^2+j}{n^2}\Big)=\sum_{j=1}^{k}\frac{j}{n^2}+\sum_{j=1}^k
\big|\Big(\log\Big(\frac{n^2+j}{n^2}\Big)-\frac{j}{n^2}\Big)\big|
\\
&=:&\frac{k^2}{2n^2}+\varepsilon'_{n,k},\nonumber
\eeq
and since $|\log(1+x)-x|\leq x^2/2$ for $x>0$ , we have that
\beqla{eq27}
|\varepsilon'_{n,k}|\leq \frac{k}{2n^2}+\sum_{j=1}^k \big| \log\big( \frac{n^2+j}{n^2}\big)-\frac{j}{n^2}\big|\leq   \frac{k}{2n^2}+\frac{1}{2n^4}\sum_{j=1}^kj^2\leq \frac{3k}{n^2}\leq \frac{6}{n}
\eeq
for $0\leq k\leq 2n.$
Since $|e^x-1|\leq 2|x|$ for $|x|\leq 1$  estimates \refeq{eq26} and \refeq{eq27} yields that
\beqla{eq28}
\lambda_{n,k}=e^{-{k^2}/{2n^2}}+ \varepsilon''_{n}\qquad {\rm with}\;\; |\varepsilon''_n|\leq 12/n.
\nonumber
\eeq

The multiplier corresponding to the sequence $(\varepsilon''_{n,k}),$ has norm less than $(2n+1)\cdot 12/n\leq 30$, since our polynomials have degree $2n.$ We next consider the main term, i.e. 
the Fourier multiplier
\beqla{eq29}
\sum_{k=0}^{\infty}c_{k}e^{ik\theta}\mapsto  \sum_{k=0}^{\infty}e^{-{k^2}/{2n^2}}c_{k}e^{ik\theta},\nonumber
\eeq
(observe that now we allow polynomials of arbitrary degree).
This map is the convolution operator $f\mapsto g* f,$ where $g$ is the positive function
\beqla{eq288}
g(\theta)=\sum_{\ell\in\integer}\sqrt{2\pi} ne^{-n^2(x-2\pi \ell )^2/2}.
\nonumber
\eeq
Especially, $\int_\torus  g (\theta )\, d\theta=1.$ The  norm of such a convolution operator  is 1 on  all the spaces $L^p(\torus ),$
and this finishes the proof of the Lemma, in view of the inequality $(30+1)/\sqrt{2\pi}\leq 20.$

Finally, the verification of  \refeq{eq21} uses the identity
\beqla{eq2222}
|P_ng((n+1)^2e^{i\theta} )|=\frac{((n+1)^2)^{(n+1)^2}}{((n+1)^2)!}\Big|
\sum_{k=1}^{2n+1}\lambda'_{n,k}b_{(n+1)^2-k}e^{-ik\theta}
\Big|, \nonumber
\eeq
with
\beqla{eq2333}
\lambda'_{n,k}:= \left(\frac{(n+1)^2}{(n+1)^2}\right)\cdots
\left(\frac{(n+1)^2-k+1}{(n+1)^2}\right). \nonumber
\eeq
The previous argument applies with minor modifications, and we obtain
\beqla{modification}
&&\| \sum_{k=1}^{2n+1}\lambda'_{n,k}b_{(n+1)^2-k}e^{-ik\theta}\|_p = \| \sum_{k=1}^{2n+1}\lambda'_{n,k}b_{(n+1)^2-k}e^{ik\theta}\|_p \leq 20
\| \sum_{k=1}^{2n+1}b_{(n+1)^2-k}e^{ik\theta}\|_p \nonumber\\
&=&20 \| \sum_{k=1}^{2n+1}b_{(n+1)^2-k}e^{-ik\theta}\|_p = 20\| \sum_{k=0}^{2n}b_{n^2+k}e^{ik\theta}\|_p\nonumber
\eeq
\end{proof}

\subsection{Proof of Theorem \ref{th:main}}We are finally ready to show that our function $f$ constructed in Section \ref{f} is has the desired properties. In what follows we show in the case $p\in (1,\infty]$ that $f$ has the desired growth and is D-frequently hypercyclic. As mentioned in the introduction, the case $p=1$ and the optimality of Theorem \ref{th:main} for any $p$ are already known after \cite{BBG-E} and \cite{Bo2}, but for the reader's convenience we present their result  using our techniques   in Remarks \ref{p=1} and \ref{optimality}.

\begin{proof}[Bounding the growth of $f$]\quad
Let first  $2\leq p\leq \infty.$ Actually,  since $\|\cdot \|_p$ increases with $p$ it is then enough to consider only the case $p=\infty .$    Assume first that $n\in {\mathcal A}_k$  is positive and even. By construction,  by Proposition \ref{heat} and by Lemma \ref{le:rudin-shapiro} we first obtain that
\beqla{bound1}
M_{P_nf}(n^2 )&\leq& 20{e^{n^2}}{n}^{-1}\| p_{\lfloor n/\alpha_k\rfloor}(z^{\alpha_k})\widetilde q_k(z)\|_\infty\leq 20\cdot 5{e^{n^2}}{n}^{-1}\sqrt{n/\alpha_k}\| \widetilde q_k\|_{\ell^1}\nonumber \\
&\leq & e^{n^2}n^{-1/2}100\ell_k\alpha_k^{-1/2}\leq c10^{-3}e^{n^2}n^{-1/2}.
\nonumber
\eeq
where one applied \refeq{imp} and the first condition in definition \refeq{alpha} of the sequence $(\alpha_k)_{k\geq 1}.$ Trivially the same bound applies for odd $n$, or for $n=0$, since then $P_nf=0.$ In a similar manner,  Proposition \ref{heat} yields that
\beqla{bound2}
M_{P_nf}((n+1)^2 )\leq c10^{-3}e^{(n+1)^2}(n+1)^{-1/2}.\nonumber
\eeq
At this stage Proposition \ref{glue} applies, and we deduce that $f$ satisfies the desired growth, i.e. $\displaystyle M_f(r)\leq ce^{r}r^{-1/4}$ for all $r>0.$

For $1<p\leq 2$ we use a similar analysis, based on replacing the polynomial $p_m$ with $p_m^*$ in (\ref{piece}), again 
with $m=\lfloor n/\alpha_k\rfloor$. Since now $p'>2$, the  above computation  takes the form
\beqla{bound3}
M_{P_nf,p}(n^2 )&\leq& 20e^{n^2}{n}^{-1}\| p^*_{\lfloor n/\alpha_k\rfloor}(z^\alpha_k)\|_p
\| \widetilde q_k\|_{\ell^1}= 20e^{n^2}{n}^{-1}\| p^*_{\lfloor n/\alpha_k\rfloor}\|_p
\| \widetilde q_k\|_{\ell^1}\nonumber \\ &\leq&20e^{n^2}{n}^{-1}\cdot 3({n/\alpha_k})^{1/p'}\| \widetilde q_k\|_{\ell^1}
\leq  60e^{n^2}n^{-1+1/p'}\ell_k\alpha_k^{-1/p'}\nonumber\\
&\leq& c10^{-3}e^{n^2}n^{-1/p}.
\nonumber
\eeq
Together with an analoguous estimate for $M_{P_nf,p}((n+1)^2 )$ we obtain the growth $\displaystyle M_{f,p}(r)\leq ce^{r}r^{-1/2p}$ again with the aid of Proposition \ref{glue}.

\medskip

\noindent{\it$f$ is D-frequently hypercyclic. }\quad This part of the argument  $f$ is independent of $p\in (1,\infty]$.  We start by observing the simple coefficient bound for the Taylor coefficients of $f:$
\beqla{size2}
|a_j|\leq j\qquad{\rm for \; all}\;\; j\geq 1.
\eeq
This is an immediate consequence of the bound \refeq{imp} together with definitions \refeq{alpha} and \refeq{piece} (resp. \refeq{piece2}) since the absolute value of the coefficients of the polynomials $p_m$ (resp. $p^*_m$) do not exceed 1, and if $j\in {\mathcal A}_k$ with $a_j\not=0$, then 
\beqla{aa}
|a_j|\leq \| \widetilde q_k\|_{\ell^1} \leq \ell_k< \alpha_k<j.\nonumber
\eeq

For any fixed even integer $n\in\mathcal A_k$ with $n\geq\alpha_k$, denote by ${\mathcal B}_n$ the set of indices $s$ such that the coefficient of $z^s$ in the polynomial $z^{n^2}p_{\lfloor 2n/\alpha_k\rfloor}(z^\alpha_k)$ is 1. With the $\{(q_k,\ell_k)\}$ used to define $f$ (so that $\ell_k\uparrow\infty$), we claim
that for each $k\geq 1$,
\beqla{near}\phantom{tttttt}
\sup_{|z|=\ell_k} | q_k(z)-\Big(\frac{d}{dz}\Big)^{s}f(z)|\leq \frac{1}{\ell_k}
\qquad {\rm for\; any}\;\;s\in {\mathcal B}_n \quad{\rm and \; for \; any}\;\, n\in {\mathcal A}_k.
\eeq
This clearly suffices 
as $\mathcal P$ dense in the space of entire functions, and obviously for any fixed $k\geq 1$ the set 
$$
\{ s\;  :\;   s\in {\mathcal B}_n,\; n\in  {\mathcal A}_k,\; n\geq\alpha_k\}
$$
has positive density since for large $n$ we have arranged that  $$\#({\mathcal B}_n)\geq (\alpha_k)^{-1}n\geq  (3\alpha_k)^{-1}\#(\{n^2,n^2+1,\ldots (n+1)^2-1\}),$$ 
and since ${\mathcal A}_k$ contains an arithmetic progression of $n$:s.

Towards \refeq{near}, fix an even $n\geq 1$ and $s\in \mathcal B_n$, where we suppose that $n\in \mathcal A_k$ with $n\geq\alpha_k.$  Our construction of $P_nf$ shows that the coefficients $a_s,a_{s+1},\ldots a_{s+ {\rm deg}\, (q_k)}$ coincide precisely with those of $q_k,$ and by the choice of $\alpha_k$ at least the next $8\ell_k$ coefficients  among the $a_j$ are zero. 
Since $(n+1)$ is odd,
we have that  $a_j=0$ for $(n+1)^2\leq j\leq (n+2)^2-1$, and so
$$
\Big(\frac{d}{dz}\Big)^{s}f(z)-q_k(z)=\sum_{j=s+8\ell_k}^{(n+1)^2-1}\frac{a_jz^{j-s}}{(j-s)!}+\sum_{j= (n+2)^2}^\infty \frac{a_jz^{j-s}}{(j-s)!} =:S_1(z)+S_2(z).
$$
By definition,  for  $n^2\leq j\leq (n+1)^2-1$ one has $|a_j|\leq \|\widetilde{q_k}\|_{\ell^1}
\leq \ell_k$,  so the first sum is bounded by
\begin{eqnarray}
\sup_{|z|=\ell_k}|S_1(z)|
\leq
\|\widetilde{q_k}\|_{\ell^1}\sum_{m=8\ell_k}^\infty\frac{\ell_k^m}{m!}
\leq 2\frac{\ell_k^{8\ell_k+1}}{(8\ell_k)!}<\frac{1}{2\ell_k}.\nonumber
\end{eqnarray}
In the above computation one observed that the ratio of any two consecutive terms in the last written series is less than 1/2,  and  the very last step was due to the estimate
\beqla{stirling}
\frac{x^{8x}}{(8x)!}\leq \frac{1}{4x^3}\qquad   {\rm for }\;\; x\geq 2.
\eeq
In turn, the inequality  \refeq{stirling}  is an easy consequence of Stirling's formula. 

In order to estimate the second sum $S_2(z)$ we write  $j=s+m$ with $m\geq (n+2)^2-(n+1)^2=2n+3$ and observe that \refeq{aa} yields that  $|a_j|=|a_{s+m}|\leq (n+1)^2+m\leq m(m-1)$. Hence we
obtain
\begin{eqnarray}
\sup_{|z|=\ell_k}|S_2(z)|&\leq &\big| \sum_{j= (n+2)^2}^\infty \frac{a_j\ell_k^{j-s}}{(j-s)!}\big| \leq
\sum_{m=2n+3}^{\infty}\frac{(n+1)^2+m}{m!}\ell_k^{m}\leq 
\ell_k^{2}\sum_{m=2n+1}^{\infty}\frac{\ell_k^{m}}{m!}\nonumber \\
&\leq &\ell_k^{2}\sum_{m=8\ell_k}^{\infty}\frac{\ell_k^{m}}{m!}\leq \frac{1}{2\ell_k}.\nonumber
\end{eqnarray}
Above one applied the knowledge $n\geq \alpha_k\geq 8\ell_k$ from  \refeq{alpha}
and the last sum  was estimated as before by  a geometric series and \refeq{stirling}. 

Put together, the estimates we obtained  for $S_1$ and $S_2$ yield \refeq{near} and the proof of Theorem \ref {th:main} is completed.
\end{proof}

\section{Remarks}\label{remarks}

\begin{remark}\label{p=1} It is instructive to analyze what happens with the above argument in the special case $p=1$. A key point in the proof is the fact that  $\widetilde{P_nf}$ contains the product $p_{[n/\alpha_k]}(z^\alpha_k)\widetilde q_k(z)$ (resp.  $p^*_{[n/\alpha_k]}(z^\alpha_k)\widetilde q_k(z)$ if $p\in (1,2)$). Hence  by choosing $\alpha_k$ large enough we may decrease the $L^p$-norm  of the first factor to compensate for the possibly increasing size of $\widetilde q_k$. However, this does not work for $p=1$ since  obviously the $L^1$-norm of any polynomial must exceed the sup-norm of its coefficients!

However, a small change in the argument produces the optimal result also in the case $p=1$ (and hence provides an alternative proof of part (iii) of Theorem \ref{th:main}). Assume that  we are given an  increasing function $\phi(r):(0,\infty)\to[1,\infty )$  with $\lim_{r\to\infty}\phi(r)=\infty$. One constructs $f$ as before in case $p\in (1,2)$, the only change is that initially in the definition of the sequence $(\alpha_k)$ one replaces condition \refeq{alpha} by the single demand $\alpha_k\geq 2d_k+8\ell_k.$ By writing  $f$ as
$$
f=\sum_{k=1}^\infty f_k\quad {\rm with}\quad f_k=\sum_{n\in {\mathcal A}_k} P_nf.
$$
the argument of Section \ref{se:proof} applies as such to each piece $f_k$: one deduces  for each $n\geq 1$ the bound $M_{P_nf_k}\leq 60e^{n^2}n^{-1}\| \widetilde q_k\|_{\ell^1},$ whence Proposition \ref{glue} yields the growth $M_{f_k,1}(r)\leq e^{|r|}|r|^{-1/2}  6\cdot 10^4\| \widetilde q_k\|_{\ell^1}$. Note that the $j$:th Taylor coefficients of $P_k$ can be nonzero only if $j\geq \alpha_k$. Hence, by further increasing  the size of $\alpha_k$,
and by recalling  \refeq{size2}  we may in addition  force $f_k$ as small as we want in any compact region. In particular, we may demand 
$$
M_{f_k,1}(r)\leq 2^{-k}\phi (r)e^{r}r^{-1/2}.
$$
Summing up we obtain the growth  $M_{f,1}(r)\leq \phi (r)e^{r}r^{-1/2}$.
\end{remark}

\begin{remark}\label{optimality}
For the reader's convenience we sketch the  proof of optimality in Theorem \ref{th:main}, although this is already contained in \cite{BBG-E}. Assume that $f$ is entire and D-frequently hypercyclic. In particular, $f$  frequently approximates the constant function  $2$ up to precision 1 in $B(0,2)$, which implies that the density of the set $H$ is positive, where $H:=\{ k\geq 1 :|a_k|\geq 1\} .$ A fortiori, there is a constant $c_1>0$ such that for infinitely many $n\geq 1$ 
\beqla{density}
\# \big( H\cap \{ n^2,n^2+1,\ldots, n^2+2n\}\big) \geq c_1n.
\eeq
Next, write down an  analogue of the identity \refeq{eq22} 
\beqla{eq123}
|f(n^2e^{i\theta} )|=\frac{(n^2)^{n^2}}{(n^2)!}\Big|
\sum_{k=-n^2}^{\infty}\lambda_{n,k}a_{n^2+k}e^{ik\theta}
\Big|, 
\eeq
where  one has extended in a natural way definition \refeq{eq23} of  $\lambda_{n,k}$ to
all values $k\geq -n^2.$ Observe that the estimates \refeq{eq26} and \refeq{eq27} show that 
$\lambda_{n,k}\geq c_2$ for $k=0,\ldots 2n$, where $c_2$ is independent of $n.$ Hence, for $p>1$ the optimality follows simply  by  considering \refeq{eq123} for arbitrarily large $n$ that satisfy \refeq{density}, applying Stirling's formula, and the standard estimate $\|\sum_{-\infty}^\infty b_ke^{ik\theta}\|_p\geq 
\|(b_k)_{k=-\infty}^{\infty}\|_{\ell^{\max (p',2)}},$ obtained from the Hausdorff-Young inequality see \cite{Ka}). In case $p=1$ one approximates similarly any given constant $A$ (instead of just $A=1$) and obtains that $\limsup_{r\to\infty}M_{f,1}(r)r^{1/2}e^{-r}\geq cA$. Since $A>1$ is arbitrary, the necessity of the increasing factor $\phi (r)$  in case (iii) of Theorem \ref{th:main} follows.

\end{remark}

\begin{remark}\label{growth}
By employing the identity \refeq{eq123} exactly as in the previous remark, and by recalling that each  set $\mathcal  A_k$ contains an arithmetic sequence of indices (it enough to consider just any single  $\mathcal  A_k$ with $q_k\not=0$), one checks that for $p\in (1,\infty]$ our function $f$ verifies the lower bound
\beqla{eq223}
M_{f,p}(r)\geq \widetilde ce^rr^{-\max (1/2p, 1/4)}\quad {\rm for}\quad r>1,
\eeq
where $\widetilde c$ is a positive constant.
Hence one has $M_{f,p}(r)\asymp e^rr^{-\max (1/2p, 1/4)}$ for all $r>1.$
However,  in general  a frequently hypercyclic function $f$ needs not  satisfy \refeq{eq223}, since one could easily modify our construction to impose infinitely many large  (dyadic) gaps to the Taylor series of $f$.
\end{remark}

\begin{remark}\label{re:reason}
A main ingredient in our construction is the selection of the right size for the blocks $P_nf$. This is actually rather delicate issue: if the blocks are suitably sparse, one obtains that the blocks $P_nf$ are independent in the sense of Proposition \ref{glue}. This can be most easily understood by comparing the growth of  a single term inside the block with the  exponential growth.  Indeed,  the ratio $r^ne^{-r}$ takes its  maximal value at $|z|=n$ and starts to decay rapidly as soon as $|z-n|>> \sqrt{n}$ (or, equivalently, $\lambda_{n,k}$ in \refeq{eq23}  decreases quickly for $k >> n$). This suggests that the  block structure should be at least as sparse as  in \refeq{blocks}, which choice also enables one to make full advantage of the use of  Lemma \ref{le:rudin-shapiro}.
On the other hand, if the   the blocks  were more sparse,  the application of Hadamard's three circles theorem in Proposition \ref{glue} would be defective.  All said,  for our proof the choice in \refeq{blocks} is essentially unique.
\end{remark}

\begin{remark}\label{hypercyclic} The question of the minimal growth of $D$-hypercyclic entire functions is considerably easier  since then there is no need to try to build in cancellation. The sharp result  \cite{S}, \cite{G-E2} says that  any $D$-frequently  hypercyclic entire function satisfies $\limsup_{r\to\infty} M_f(r)e^{-r}\sqrt{r} =\infty$ and this result is optimal.
\end{remark}


\end{document}